\documentclass[11pt, a4paper]{article}
\usepackage{latexsym}
\usepackage{indentfirst}
\usepackage{graphicx}
\usepackage{amsmath, amsthm, amsfonts, amssymb}	
\usepackage{hyperref}

\newtheorem{lem}{Lemma}[section]

\newtheorem{thm}{Theorem}[section]

\numberwithin{equation}{section}

\begin{document}
\title{A goodness-of-fit test of the errors in nonlinear autoregressive time series models with stationary $\alpha$-mixing error terms}
\author{Kyong-Hui Kim, Myong-Guk Sin, Ok-Kyong Kim \\
Faculty of Mathematics, \textbf{Kim Il Sung} University, \\
Pyongyang, Democratic People's Republic of Korea\\
kim.kyonghui@yahoo.com
}

\date{}
\maketitle

\begin{abstract}
In this work we deal with the problem of fitting an error density to the goodness-of-fit test of the errors in nonlinear 
autoregressive time series models with stationary $\alpha$-mixing error terms. The test statistic is based 
on the integrated squared error of the nonparametric error density estimate and the null error density. 
By deriving the asymptotic normality of test statistics in these models, we extend the result of Cheng and Sun 
(Statist. Probab. Lett. \textbf{78}, 1(2008), 50-59) in the model with i.i.d error terms to the more general case.
\end{abstract}

%
%

\section{Introduction}
The purpose of the goodness-of-fit tests is to test hypotheses on the empirical distributions fitting some theoretical law. 
Much recent work has been devoted to the goodness-of-fit tests of the errors in variant models; see for example 
\cite{amj, che, kla}. 
In autoregressive time series models, the goodness-of-fit tests based on the residual empirical process 
have been extensively studied; for more details concerning them, we refer to \cite{bol, kou}.
Bachmann and Dette \cite{bac} studied the Bickel-Rosenblatt test by considering the asymptotic behaviour of the test 
statistic under a fixed alternative.  
They proved that, under such conditions, a standardized version of the Bickel-Rosenblatt test statistic based on i.i.d. 
observations is asymptotically normal distributed, but with a different rate of convergence.
Cheng and Sun \cite{che} derived asymptotic normality of the Bickel-Rosenblatt test statistic in nonlinear 
autoregressive time series models with i.i.d. errors. 
We extend this result to nonlinear autoregressive time series models with stationary $\alpha$-mixing error terms.

Let $\{X_i, i=0, \pm 1, \pm 2, \cdots \}$ be a strictly stationary process of real random variables 
obeying the model
\begin{equation} \label{eq1.1}
X_i=r_{\theta}(X_{i-1}, \cdots, X_{i-p})+\varepsilon_i 
\end{equation}
  for some $\theta=(\theta_1, \cdots, \theta_q)^T\in\it\Theta\subset R^q$, where 
($r_{\theta}; \theta \in \it \Theta$) is a family of known measurable functions from $R^p$ to $R$. 
Unlike in \cite{che}, we assume that the errors $\varepsilon_i$ are $\alpha$-mixing random variables with common 
density $f$ and $X_{i-1}, \cdots, X_{i-p}$ are independent of $\{ \varepsilon_i, i=1, 2, \cdots \}$. 
We focus on the problem of testing the hypothesis 	
\begin{equation} \label{eq1.2}
H_0 : f=f_0 \quad \textnormal{vs.} \quad H_1 : f \neq f_0, 
\end{equation}
where $f_0$ is a prescribed density based on the data $\{ X_{1-p}, \cdots, X_0, X_1, \cdots, X_n \}$.

We perform a test using the integrated square deviation of a kernel type density estimator 
based on the residuals from the expectation of the kernel error density based on the true errors. 
Let $\hat{\theta}=(\hat{\theta}_1, \cdots, \hat{\theta}_q)^T$ be an estimator of $\theta$, 
and define the residuals for $i=1, 2, \cdots$,
\begin{equation*}
\hat{\varepsilon}_i=X_i-r_{\hat{\theta}}(X_{i-1}, \cdots, X_{i-p}).
\end{equation*}
For a kernel density function $K$, the kernel type estimator of the error density $f(t)$ is defined as 
\begin{equation*}
\hat{f}_n(t) = \frac{1}{n} \sum_{i=1}^n K_{h_n}(t-\hat{\varepsilon}_i), t \in R,
\end{equation*}
 where $K_h(\cdot)=(1/h)K(\frac{\cdot}{h})$ is a scaled kernel and $h=h_n$ denotes a bandwidth tending to zero. 
We also define the kernel error density based on the true errors $\varepsilon_1, \cdots, \varepsilon_n$, 
which we cannot observe, as follows:
 \begin{equation*}
f_n(t) = \frac{1}{n} \sum_{i=1}^n K_{h_n}(t-\varepsilon_i), t \in R.
\end{equation*}

For the problem of testing the hypothesis (1.2) we use the integrated squared deviation of $\hat{f}_n$ from 
\begin{equation*}
Ef_n(t) = \int K(x)f(t-h_nx)dx = K_h \ast f(t),
\end{equation*}
where $K_h \ast f$ denotes the convolution of the functions $K_h$ and $f$, i.e., we reject the null-hypothesis 
$H_0 : f=f_0$ for large values of the statistic
\begin{equation*}
\hat{T}_n=\int [\hat{f}_n(t)-K_h \ast f_0(t)]^2 dt.
\end{equation*}
This $\hat{T}_n$ is an analogue of the Bickel-Rosenblatt statistic proposed in the case of the 
observable $\varepsilon_i$'s 
\begin{equation*}
T_n=\int [f_n(t)-E(f_n(t))]^2 dt,
\end{equation*}
see \cite{bic}.

%
%

\section{Basic assumptions and preliminaries}

In this section we introduce some basic assumptions on the nonlinear autoregressive model \eqref{eq1.1} and the estimator 
and give some preliminaries which can be used to prove our main results. The same assumptions on the autoregression 
function $r_\theta$ and the estimator $\hat{\theta}$ for $\theta$ as in \cite{che} are adopted here. 
Throughout the paper we assume that limits are taken as $n \to \infty$ unless otherwise specified.\\

$(\textnormal{A}_1)$. Let $U \subset \Theta \subset R^q$ be an open neighborhood of $\theta$. 
We assume that, for all $y \in \mathbb{R}^p, \vartheta = (\vartheta_1, \cdots, \vartheta_q) \in U$ and $j, k=1, \cdots, q$,
\begin{equation*}
\left| \frac{\partial}{\partial \vartheta_j}r_{\vartheta}(y) \right| \leq M_1(y)
\end{equation*}
\begin{equation*}
\left| \frac{\partial^2}{\partial \vartheta_j \partial \vartheta_k} r_{\vartheta}(y) \right| \leq M_2(y),
\end{equation*}     
where $EM_1^4(X_{i-1}, \cdots, X_{i-p})<+\infty$ and $EM_2^4(X_{i-1}, \cdots, X_{i-p})<+\infty$ 
for  $i \geq 1$. \\

For all $1 \leq i \leq n$  and $1 \leq j \leq q$, let 
\begin{equation*}
Y_{ij}= \frac{\partial}{\partial \theta_j} r_{\theta}(X_{i-1}, \cdots, X_{i-p}).\\
\end{equation*}     

$(\textnormal{A}_2)$. We assume that there exists $\alpha<1$  such that $Y_{ij}$'s satisfy
\begin{equation*}
\sum_{i=1}^n Y_{ij} = \textnormal{O}_\textnormal{p}(n^{\alpha}), ~ j=1, 2, \cdots, q. \\
\end{equation*}    
                      
$(\textnormal{A}_3)$. We assume that the estimator $\hat{\theta} = (\hat{\theta}_1, \cdots, \hat{\theta}_q)^T$ 
for $\theta$ (based on $X_0, X_1, \cdots, X_n$) satisfies the law of iterated logarithm, i.e., there exists 
a constant $C_1 (0<C_1<\infty)$ such that 
\begin{equation*}
\lim \sup_{n \to \infty} \sqrt{\frac{n}{\log(\log n)}} \left \vert \hat{\theta}-\theta \right \vert \leq C_1,
\end{equation*}    
 where $\vert \hat{\theta}-\theta \vert = \sqrt{\sum_{j=1}^q (\hat{\theta}_j-\theta_j)^2}$. \\

In this work we derive the asymptotic distribution of $\hat{T_n}$  under $H_0$. 
In order to calculate the probability of type II 
error when $\hat{T_n}$  is used to test hypothesis  \eqref{eq1.2}, we consider the asymptotic distribution of $\hat{T_n}$ 
under one fixed alternative in $H_1$ of  \eqref{eq1.2} in the sense of
\begin{equation*}
d(f, f_0)=\int (f-f_0)^2(x)dx>0.
\end{equation*}    

Next we describe some basic assumptions on the error density  $f$, the kernel density $K$ and the bandwidth $h_n$. 
	
(D). $f$  is two time continuously differentiable with bounded first and second derivatives, and $f^2$  is integrable. 

(K). $(\textnormal{K}_1) ~ K$ is a continuous bounded symmetric kernel with compact support.

\quad \quad $(\textnormal{K}_2)  ~K'''$ exists and is bounded. $K', (K')^2, K''$ and $(K'')^2$ are integrable.

(H).  $nh_n^2 \to \infty $ and $h_n \to 0$.\\

Note that assumption $(\textnormal{K}_1)$ implies that 
\begin{equation*}
\int x^2K(x)dx<\infty \quad \textnormal{and} \quad \int K^2(x)dx<\infty.\\
\end{equation*}    

Under the above assumptions, Cheng and Sun \cite{che} established the following results.
	
\begin{lem} \label{lem2.1}
Under assumptions $(\textnormal{A}_1)$ and $(\textnormal{A}_3)$, we have
\begin{equation*}
\sum_{i=1}^n(\hat{\varepsilon}_i-\varepsilon_i)^2 = \textnormal{O}_\textnormal{p} ( \textnormal{log(log} n)).
\end{equation*}
\end{lem}
 
\begin{lem} \label{lem2.2}
Under assumptions \textnormal{(D)} and \textnormal{(K)}, we have
\begin{eqnarray*}
& \textnormal{(i)} & \int \Big[ E \Big( K'\Big( \frac{t-\varepsilon_i}{h_n}\Big) \Big) \Big]^2dt =  \textnormal{O}(h_n^2), \quad 
\int E \Big( K'\Big( \frac{t-\varepsilon_1}{h_n}\Big) \Big)^2dt =  \textnormal{O}(h_n) \\
& \textnormal{(ii)} & \int \Big[ E \Big( K''\Big( \frac{t-\varepsilon_i}{h_n}\Big) \Big) \Big]^2dt =  \textnormal{O}(h_n^2), \quad 
\int E \Big( K''\Big( \frac{t-\varepsilon_1}{h_n}\Big) \Big)^2dt =  \textnormal{O}(h_n). 
\end{eqnarray*}
\end{lem}
 
\begin{lem} \label{lem2.3}
Suppose that assumptions $(\textnormal{A}_1)-(\textnormal{A}_3), (\textnormal{D}), (\textnormal{K})$ 
and $(\textnormal{H})$ hold and the bandwidth $h_n$ satisfies the following condition
\begin{equation*}
n^{-1}h_n^{-4}\textnormal{(log(log} n))^2 \to 0
\end{equation*}
and, moreover,
\begin{equation*}
n^{2(\alpha-1)}h_n^{-3/2} \textnormal{log(log} n) \to 0.
\end{equation*}
Then we have
\begin{equation} \label{eq2.1}
\int \Big[ \hat{f}_n(t)-f_n(t) \Big]^2 dt = \textnormal{O}_ \textnormal{p} \Big( \frac{(\textnormal{log(log} n))^2}{n^2h_n^4}+\frac{\textnormal{log(log} n)}{n^{3-2\alpha}h_n^2} \Big) = \textnormal{o}_ \textnormal{p} \Big(\frac{1}{n \sqrt{h_n}} \Big)
\end{equation}
\end{lem}

%
%
	
\section{Main results}

In this section we derive the asymptotic normality of the Bickel-Rosenblatt test statistic in nonlinear 
autoregressive time series models with stationary $\alpha$-mixing error terms.\\

We start with the following property of stationary $\alpha$-mixing random variables. 
\begin{lem} \label{lem3.1}
Suppose that the stationary sequence $\{ X_i \}$ satisfies $\alpha$-mixing condition. If the random variables $\xi$ and $\eta$ 
are measurable for $\Im_{-\infty}^t$ and $\Im_{t+\tau}^{\infty}$ and $\vert \xi \vert <C_1, \vert \eta \vert <C_2$  
then we obtain 
\begin{equation*}
\vert E\xi\eta-E\xi \cdot E\eta \vert \leq 4C_1C_2 \alpha(\tau).
\end{equation*}
\end{lem}
The proof of this lemma is simple, so is omitted.\\
 
We are now in position to formulate our main results in this exposition.
%
%
\begin{thm} \label{thm3.1}
Suppose that assumptions \textnormal{(D), (K)} and \textnormal{(H)} are satisfied. Then Bickel-Rosenblatt test statistics 
\begin{equation*}
T_n=\int \big[ f_n(t)-K_h \ast f_0(t) \big]^2 dt
\end{equation*}    
has the following properties:

\textnormal{(i)} Under the null hypothesis $H_0 : f=f_0$, as $n \to \infty$  
\begin{equation} \label{eq3.1}
n \sqrt h_n \Big[ T_n - \frac{1}{nh_n} \int K^2(x)dx \Big] \to N\Big(0, 2 \int f_0^2(x)dx \int (K \ast K)^2(x)dx \Big).
\end{equation}    

\textnormal{(ii)} Under the alternative $H_1 : f \neq f_0$, as $n \to \infty$  
\begin{equation} \label{eq3.2}
\sqrt n \Big[ T_n -\int (K_h\ast (f-f_0))^2(x)dx \Big] \to N\Big(0, 4 Var[(f-f_0)(\varepsilon_1)] \Big).
\end{equation}    
\end{thm}
	
\textbf{Proof} Recall that we are establishing the asymptotic normality under the null hypothesis $f=f_0$ 
and under fixed alternatives $f \neq f_0$ with different rates of convergence in both cases.
Let $f$ denote the \textquotedblleft true\textquotedblright  density of the random variables $\varepsilon_i$. By the definition of the statistic $T_n$ 
and the density estimate $f_n$, we obtain the following decomposition:
\begin{align*}
&T_n  =  \int [f_n-K_h \ast f_0]^2(x) dx  \nonumber\\
	& \quad =  \int [f_n-K_h \ast f]^2(x) dx +2 \int [f_n-K_h \ast f](x) g_h(x)dx + \int g_h^2(x)dx \nonumber \\
	&  \quad =  \frac{2}{n^2} \sum_{i<j} \int [K_h(x-\varepsilon_i)-e_h(x)][K_h(x-\varepsilon_i)-e_h(x)]dx \nonumber\\
	&   \quad  \quad +\frac{2}{n} \sum_{i=1}^n [(K_h \ast g_h)(\varepsilon_i)-E[(K_h \ast g_h)(\varepsilon_i)]] \nonumber\\
	&   \quad  \quad +\frac{1}{n^2} \sum_{i=1}^n [K_h (x-\varepsilon_i)-e_h(x)]^2 dx + g_h^2(x)dx.
\end{align*}
where the functions $e_h, g_h$ are defined by $e_h=K_h \ast f$ and $g_h=K_h \ast (f-f_0)$, respectively. 
Simple calculation implies
\begin{equation*}
\frac{1}{n^2} \sum_{i=1}^n \int \big[ K_h(x-\varepsilon_i)-e_h(x) \big]^2 dx = 
\frac{1}{nh} \int K^2(x)dx+\textnormal{O}_ \textnormal{p} \Big( \frac{1}{n} \Big),
\end{equation*}
 and therefore we have the stochastic expansion
\begin{eqnarray*}
T_n & - & \frac{1}{nh} \int K^2(x)dx - \int [K_h \ast (f-f_0)]^2(x)dx \nonumber\\
	& = & \frac{2}{n^2} \sum_{i<j} H_n(\varepsilon_i, \varepsilon_j)+\frac{2}{n} \sum_{i=1}^n Y_i+ 
\textnormal{O}_ \textnormal{p} \Big( \frac{1}{n} \Big),
\end{eqnarray*}
where 
\begin{eqnarray*}
H_n(\varepsilon_i, \varepsilon_j) & = & \int \big[K_h(x-\varepsilon_i)-e_h(x) \big] \big[K_h(x-\varepsilon_j)-e_h(x) \big]dx, \\
Y_i & = & (K_h \ast g_h)(\varepsilon_i)-E \big[K_h \ast g_h(\varepsilon_i) \big].
\end{eqnarray*}
      	
Denote the first term in this decomposition as 
\begin{equation*}
U_n=\frac{2}{n^2} \sum_{i<j} H_n(\varepsilon_i, \varepsilon_j).
\end{equation*}
Note that $U_n$ does not depend on the density $f_0$ specified by the null hypothesis. 
It is clear that $H_n$ is symmetric and
\begin{equation*}
\lim_{n \to \infty} E\big[ H_n(\varepsilon_1, \varepsilon_2) \vert \varepsilon_1 \big]=0, \quad
\lim_{n \to \infty} E\big[ H_n^2(\varepsilon_1, \varepsilon_2)\big]< \infty
\end{equation*}
for each $n \in \mathbb{N}$. In fact,
\begin{equation*}
E\big[ H_n(\varepsilon_i, \varepsilon_j) \vert \varepsilon_i \big] = \int E \big[ (K_h(x- \varepsilon_i)-e_h(x)) (K_h(x- \varepsilon_i)-e_h(x)) \vert \varepsilon_i \big]dx.
\end{equation*}
Denote the random variable in the integrate symbol by $\eta$, then we have 
\begin{equation*}
E \left \vert E \big\{\eta \vert \Im_{-\infty}^0 \big\} -E\eta \right \vert = E \big\{ \xi_1 \big(E\big(\eta \vert \Im_{-\infty}^0 \big)-E\eta \big) \big\},
\end{equation*}
where $\xi_1= \textnormal{sgn} \big(E\big(\eta \vert \Im_{-\infty}^0 \big)-E\eta \big)$ and it is measurable for 
$\Im_{-\infty}^0$. It follows (via $|\eta| \leq 4$) that
\begin{equation*}
\left \vert E \xi_1 \eta-E \xi_1E\eta \right\vert \leq 4 \left \vert E \xi_1 \eta_1-E \xi_1E\eta_1 \right\vert,
\end{equation*}
where $\eta_1= \textnormal{sgn} \big(E\big\{\xi_1 \vert \Im_{-\infty}^0 \big\}-E\xi_1 \big)$.

By Lemma \ref{lem3.1}, we have 	
\begin{eqnarray*}
\left \vert E \xi_1 \eta_1-E \xi_1E\eta_1 \right\vert \leq 4 \alpha(\tau),
\end{eqnarray*}
and therefore
\begin{equation*}
E \left \vert E \big\{\eta \vert \Im_{-\infty}^0 \big\} -E\eta \right \vert \leq 16 \alpha(\tau),
\end{equation*}
where $\tau = |i-j|$.
Since $\{\varepsilon_i\}$ is a strictly mixing sequence with coefficient $\alpha(\tau)$, the left-hand side of above equation 
converges zero as $\tau \to \infty$. And we have
\begin{equation*}
EK_h(x-\varepsilon_i) = K_h \ast f(x) = e_h.
\end{equation*}
Let 
\begin{equation*}
\xi = K_h(x-\varepsilon_i), \eta = K_h(x-\varepsilon_j),
\end{equation*}
then by Lemma \ref{lem3.1} we obtain  $|E \eta| \leq 4 \alpha(\tau)$, which implies that 	
\begin{equation*}
\lim_{n \to \infty} E[H_n(\varepsilon_i, \varepsilon_j) \vert \varepsilon_j] = 0.
\end{equation*}
The other moment limit results can be proved in the same way. Applying the central limit theorem for degenerate 
U-statistics completes the proof. $\Box$ 
	
%
%

\begin{thm}  \label{thm3.2}
Suppose that assumptions $(\textnormal{A}_1)-(\textnormal{A}_3), (\textnormal{D}), (\textnormal{K})$ 
and $(\textnormal{H})$ are satisfied and that the bandwidth $h_n$ satisfies the following:
\begin{equation} \label{eq3.3}
n^{2(\alpha-1)}h_n^{-2}\textnormal{log(log} n) \to 0,
\end{equation}    
\begin{equation} \label{eq3.4}
n^{-1}h_n^{-4}\textnormal{(log(log} n))^2 \to 0.
\end{equation}    
Then the test statistics $\hat{T}_n$ has the following properties:

\textnormal{(i)} Under the null hypothesis $H_0 : f=f_0$, as $n \to \infty$,  
\begin{equation*}
n \sqrt {h_n} \Big[ \hat{T}_n - \frac{1}{nh_n} \int K^2(x)dx \Big] \to N\Big(0, 2 \int f_0^2(x)dx \int (K \ast K)^2(x)dx \Big).
\end{equation*}    

\textnormal{(ii)} Under the alternative $H_1 : f \neq f_0$, as $n \to \infty$,  
\begin{equation*}
\sqrt {n} \Big[ \hat{T}_n -\int (K_h\ast (f-f_0))^2(x)dx \Big] \to N\Big(0, 4 Var[(f-f_0)(\varepsilon_1)] \Big).
\end{equation*}    
\end{thm}

\textbf{Proof} First we prove (i). By \eqref{eq3.1}, it sufficies to show that 
\begin{equation*}
n \sqrt{h_n} (\hat{T}_n-T_n) = \textnormal{o}_{\textnormal{p}}(1)
\end{equation*}
 From the definition of $\hat{T}_n$ and $T_n$, we obtain
\begin{equation} \label{eq3.5}
\left \vert \hat{T}_n-T_n \right \vert \leq \int \big( \hat{f}_n(t)-f_n(t) \big)^2dt + 2 \Big[ \int \big( \hat{f}_n(t)-f_n(t) \big)^2dt \Big]^{1/2} \sqrt{T_n}.
\end{equation}
Using Lemma \ref{lem2.2} and the fact
\begin{equation*}
T_n = \textnormal{O}_{\textnormal{p}} \Big( \frac{1}{nh_n} \Big)
\end{equation*}
obtained from  \eqref{eq3.1}, we have
\begin{eqnarray*}
\left \vert \hat{T}_n-T_n \right \vert & = & \textnormal{o}_{\textnormal{p}} \Big( \frac{1}{n \sqrt{h_n}} \Big) + 
\textnormal{O}_{\textnormal{p}} \left( \frac{1}{n\sqrt{h_n}} \sqrt{\frac{(\textnormal{log(log} n))^2}{nh_n^4} 
+ \frac{\textnormal{log(log} n)}{n^{2-2\alpha}h_n^2}} \right) \\
	& = & \textnormal{o}_{\textnormal{p}} \Big( \frac{1}{n \sqrt{h_n}} \Big).
\end{eqnarray*}
Here \eqref{eq3.3} and \eqref{eq3.4} were also used. This completes the proof of (i).

Next we prove (ii). By \eqref{eq3.2}, it suffices to show that
\begin{equation*}
\sqrt{n}(\hat{T}_n-T_n) = \textnormal{o}_{\textnormal{p}}(1).
\end{equation*}
Again by \eqref{eq3.2} we obtain
\begin{equation} \label{eq3.6}
T_n = \textnormal{O}_{\textnormal{p}}(1),
\end{equation}
and from \eqref{eq3.5}, \eqref{eq3.6} and \eqref{eq2.1}, it follows that
\begin{eqnarray*}
\left \vert \hat{T}_n-T_n \right \vert & = & \textnormal{o}_{\textnormal{p}} \Big( \frac{1}{\sqrt{n} \sqrt{nh_n}} \Big) + 
\textnormal{O}_{\textnormal{p}} \left( \frac{1}{\sqrt{n}} \sqrt{\frac{(\textnormal{log(log} n))^2}{nh_n^4} 
+ \frac{\textnormal{log(log} n)}{n^{2-2\alpha}h_n^2}} \right) \\
	& = & \textnormal{o}_{\textnormal{P}} \Big( \frac{1}{\sqrt{n}} \Big).
\end{eqnarray*}
Here we also used \eqref{eq3.3}, \eqref{eq3.4} and the fact $nh_n \to \infty$, which is guaranteed by \eqref{eq3.4} and assumption $H$. 
Now the proof of (ii) is straightforward. $\Box$


\begin{thebibliography}{99}

\bibitem{amj}
A. Amjad, Salahuddin, Alamgir, \textit{Testing goodness-of-fit in autoregressive fractionally integrated moving-average models 
with conditional hetroscedastic errors of unknown form}, Res. J. Recent Sci., \textbf{2}, 5(2013), 39-43.

\bibitem{bac}
D. Bachmann, H. Dette, \textit{A note on the Bickel-Rosenblatt test in autoregressive time series}, 
Statist. Probab. Lett. \textbf{74}, 3(2005), 221-234.

\bibitem{bic}
P. J. Bickel, M. Rosenblatt, \textit{On some global measures of the deviations of density function estimators}, 
Ann. Statist. \textbf{1}, 6(1973), 1071-1095.
	
\bibitem{bol}
M. V. Boldin, P. I. Nelson, \textit{On conditional least squares estimation for stochastic processes}, 
Ann. Statist. \textbf{6}, 3(1978), 629-642.

\bibitem{che}
F. Cheng, S. Sun, \textit{A goodness-of-fit test of the errors in nonlinear autoregressive time series models}, 
Statist. Probab. Lett. \textbf{78}, 1(2008), 50-59.

\bibitem{kla}
B. Klar, F. Lindner, S.G. Meintanis, \textit{Specification tests for the error distribution in GARCH models}, Comput. Statist. Data Anal., 
\textbf{56}, 11(2012), 3587-3598.

\bibitem{kou}
H. L. Koul, \textit{Weighted empirical processes in dynamic nonlinear models}, 
Lecture Notes in Statistics, vol.166, Springer, New York, 2002.

\end{thebibliography}
\end{document}